\documentclass[reqno]{amsart}

\usepackage{amsmath}
\usepackage{amssymb}
\usepackage{amscd}

\title{Toposes from forcing for intuitionistic ZF with atoms}


\author{Keita Yamamoto}

\email{yamamoto.nnath@gmail.com}

\thanks{During the preparation of this paper I was supported by my father.}

\newtheorem{theorem}{Theorem}[section]
\newtheorem{proposition}[theorem]{Proposition}

\newtheorem{corollary}[theorem]{Corollary}
\theoremstyle{definition}
\newtheorem{definition}[theorem]{Definition}

\newcommand{\setbegin}{\left\{\,}
\newcommand{\setmid}{\mathrel{}\middle|\mathrel{}}
\newcommand{\setend}{\,\right\}}

\DeclareMathOperator{\Hom}{Hom}
\DeclareMathOperator{\dom}{dom}
\DeclareMathOperator{\cod}{cod}
\DeclareMathOperator{\ama}{ama}
\DeclareMathOperator{\Ob}{Ob}
\DeclareMathOperator{\Arr}{Arr}
\DeclareMathOperator{\Sh}{Sh}

\usepackage{hyperref}

\begin{document}
\begin{abstract}
We introduce the forcing model of IZFA (Intuitionistic Zermelo-Fraenkel set theory with Atoms) for every Grothendieck topology and prove that the topos of sheaves on every site is equivalent to the category of `sets in this forcing model'.
\end{abstract}
\maketitle

\section{Introduction}
For a complete Heyting algebra $H$, the \textit{Heyting-valued model} $V^{(H)}$ \textit{of Intuitionistic Zermelo-Fraenkel set theory} (IZF) is obtained by carrying out the definition of the Boolean-valued model $V^{(B)}$ of ZFC with $H$ in place of a complete Boolean algebra $B$.
Then it can be shown \cite[pp. 179--181]{bell} that the topos $\Sh (H)$ of sheaves on $H$ is equivalent to the category $\mathrm{Set}^{(H)}$ of `sets in $V^{(H)}$',
which is defined more precisely as follows:
\begin{itemize}
\item we identify elements $u, v$ of $V^{(H)}$ when the truth value $\| u = v \|_{V^{(H)}} \in H$ is equal to 1,
\item the objects of $\mathrm{Set}^{(H)}$ are the (identified) elements of $V^{(H)}$,
\item the arrows of $\mathrm{Set}^{(H)}$ are those (identified) elements $f$ of $V^{(H)}$ for which $\| f \text{ is a function} \|_{V^{(H)}}$ = 1.
\end{itemize}
In this paper, for every Grothendieck topology $J$ on every small category $\mathcal{C}$, we introduce the \textit{forcing model of} IZFA (\textit{Intuitionistic Zermelo-Fraenkel set theory with Atoms}) as an extended version of Heyting-valued models of IZF and prove that the topos $\Sh (\mathcal{C}, J)$ of sheaves on $(\mathcal{C}, J)$ is equivalent to the category $\mathrm{Set}^{(\mathcal{C}, J)}$ of `sets in this forcing model'.

This forcing for IZFA is a modification of forcing for IZF in \cite{streicher}.
The points of modification are as follows:
\begin{enumerate}
\item the universe of the forcing model of IZFA includes the arrows of $\mathcal{C}$ as atoms while that of \cite{streicher} is without atoms,
\item the forcing model of IZFA is defined without using toposes directly since this formulation is more convenient for the author and for other set theorists than that of \cite{streicher}.
\end{enumerate}
The point (1) is necessary to prove that the categories $\Sh (\mathcal{C}, J)$ and $\mathrm{Set}^{(\mathcal{C}, J)}$ are equivalent, which is the main theorem (Theorem \ref{thm:equiv}).
The point (2) will enable category theorists and set theorists to communicate more with each other.

As a related work, it is shown in \cite{awodey} that every Grothendieck topos has a equivalent topos which is the universe of some model of IZFA.
Our result is stronger than it since our forcing model has only set many atoms while the model in \cite{awodey} has proper class many atoms.

In section \ref{sec:forcing}, we define forcing for IZFA and present some propositions on it.
In section \ref{sec:toposes}, we define the category $\mathrm{Set}^{(\mathcal{C}, J)}$ for each site $(\mathcal{C}, J)$ and prove the main theorem.

Notation and terminology:
\begin{itemize}
\item On Grothendieck topologies or sheaves, we adopt the terminology of \cite[Chapter III]{sgl}.
\item $\Ob (\mathcal{C})$ is the class of all objects of a category $\mathcal{C}$.
\item $\Arr (\mathcal{C})$ is the class of all arrows of a category $\mathcal{C}$.
\item $\Hom_{\mathcal{C}}(\mathrm{any}, B) := \bigcup_{A \in \Ob (\mathcal{C})} \Hom_{\mathcal{C}} (A, B)$.
\item $\Hom_{\mathcal{C}}(A, \mathrm{any}) := \bigcup_{B \in \Ob (\mathcal{C})} \Hom_{\mathcal{C}} (A, B)$.
\item $\mathcal{L}_{\in}$ is the first-order language with two binary predicate symbols $=$ (equality), $\in$ (membership).
\item $\mathcal{L}_{\mathrm{atom}}$ is the first-order language obtained by adding two unary predicate symbols $*\colon\mathrm{atom}$, $*\colon\mathrm{set}$ to $\mathcal{L}_{\in}$.
\end{itemize}

\section{Forcing for IZFA}\label{sec:forcing}
Let $(\mathcal{C}, J)$ be a site.\\
In this section, we introduce the forcing model $\left(W^{(\mathcal{C}, J)}, \Vdash_{(\mathcal{C}, J)}\right)$, which consists of the class-valued presheaf $W^{(\mathcal{C}, J)}$ and the forcing relation $\Vdash_{(\mathcal{C}, J)}$.
The definition of this forcing is a modification of forcing for IZF in \cite{streicher}.
After giving the definition, we present some propositions on it, which are used in the next section.
Most proofs of these propositions are omitted in this paper since we can prove them almost by arguments similar to that of forcing for ZFC with posets familiar to set theorists.

\subsection{Definition of forcing}
We fix two injective class functions $x \mapsto x^{(\mathrm{atom})}$ and $x \mapsto x^{(\mathrm{set})}$ on $V$ whose ranges $\left\{ x^{(\mathrm{atom})} \setmid x \in V  \right\}$ and $\left\{ x^{(\mathrm{set})} \setmid x \in V  \right\}$ are disjoint.\\
For example, $x^{(\mathrm{atom})} := (x, 0)$ and $x^{(\mathrm{set})} := (x, 1)$.

\begin{definition}
We define a presheaf $W_{\alpha}^{(\mathcal{C}, J)} \colon \mathcal{C}^{\mathrm{op}} \rightarrow \mathrm{Set}$ for each ordinal $\alpha$ by transfinite recursion as follows: 

[ Case: $\alpha = 0$ ]
For $A \in \mathrm{Ob}(\mathcal{C})$,
\[
W_0^{(\mathcal{C}, J)} (A) := \left\{ k^{(\mathrm{atom})} \setmid k \in \Hom_{\mathcal{C}} \left( A, \mathrm{any} \right)  \right\}.
\]
For $f \in \Hom_{\mathcal{C}} ( A, B )$, we define a function $W_0^{(\mathcal{C}, J)} (f) \colon W_0^{(\mathcal{C}, J)} (B) \rightarrow W_0^{(\mathcal{C}, J)} (A)$ by
\[
W_0^{(\mathcal{C}, J)} (f) \left( k^{(\mathrm{atom})} \right) := \left( k \circ f \right)^{(\mathrm{atom})}.
\]

[ Case: successor ordinal $\alpha + 1$ ]
For $A \in \mathrm{Ob}(\mathcal{C})$, we define $\widetilde{W}_{\alpha + 1}^{(\mathcal{C}, J)} (A)$ to be the set of all $a^{(\mathrm{set})}$ satisfying two conditions (1) and (2):
\begin{enumerate}
\item $a \subseteq \bigcup_{f \in \Hom_{\mathcal{C}} (\mathrm{any}, A)} W_\alpha^{(\mathcal{C}, J)} ( \dom f ) \times \{ f \}$,
\item $\left( W_\alpha^{(\mathcal{C}, J)} (g)(b), f \circ g \right) \in a$ for all $(b, f) \in a$ and all $g \in \Hom_\mathcal{C} ( \mathrm{any}, \dom f )$.
\end{enumerate}
Let $W_{\alpha + 1}^{(\mathcal{C}, J)} (A) := \widetilde{W}_{\alpha + 1}^{(\mathcal{C}, J)} (A) \cup W_0^{(\mathcal{C}, J)} (A)$. 
For $f \in \Hom_{\mathcal{C}} ( A, B )$, we define a function $W_{\alpha + 1}^{(\mathcal{C}, J)} (f) \colon W_{\alpha + 1}^{(\mathcal{C}, J)} (B) \rightarrow W_{\alpha + 1}^{(\mathcal{C}, J)} (A)$ by:
\begin{itemize}
\item $W_{\alpha + 1}^{(\mathcal{C}, J)} (f) \left( x^{(\mathrm{set})} \right) := \left\{ (y, g) \setmid g \in \Hom_\mathcal{C} (\mathrm{any}, A), (y, f \circ g) \in x \right\}^{(\mathrm{set})}$,
\item $W_{\alpha + 1}^{(\mathcal{C}, J)} (f) \left( k^{(\mathrm{atom})} \right) := W_0^{(\mathcal{C}, J)} (f) \left( k^{(\mathrm{atom})} \right)$.
\end{itemize}

[ Case: limit ordinal $\gamma$ ]
For $A \in \mathrm{Ob}(\mathcal{C})$,
\[
W_\gamma^{(\mathcal{C}, J)} (A) := \bigcup_{\alpha < \gamma} W_\alpha^{(\mathcal{C}, J)} (A).
\]
For $f \in \Hom_{\mathcal{C}} ( A, B )$, since the functions $\left\{ W_{\alpha}^{(\mathcal{C}, J)} (f) \setmid \alpha < \gamma \right\}$ are pairwise compatible by the definition, we define
\[
W_\gamma^{(\mathcal{C}, J)} (f) := \bigcup_{\alpha < \gamma} W_\alpha^{(\mathcal{C}, J)} (f).
\]
\end{definition}

\begin{definition}
We define $W^{(\mathcal{C}, J)}(A) := \bigcup_{\alpha \in \mathrm{Ord}} W_\alpha^{(\mathcal{C}, J)} (A) $ for $A \in \mathrm{Ob}(\mathcal{C})$.
We will use dotted letters $\dot{a}, \dot{b}, \dot{c} \dots$ to denote elements of $W^{(\mathcal{C}, J)}(A)$.
For $\dot{a} \in W^{(\mathcal{C}, J)}(A)$ and $f \in \Hom_\mathcal{C} (\mathrm{any}, A)$, we define $\dot{a} \cdot f$ to be $W_\alpha^{(\mathcal{C}, J)} (f) (\dot{a})$ for some ordinal $\alpha$ for which $\dot{a} \in W_\alpha^{(\mathcal{C}, J)} (A) $.
This definition is independent of choice of $\alpha$ since the functions $\left\{ W_{\alpha}^{(\mathcal{C}, J)} (f) \setmid \alpha \in \mathrm{Ord} \right\}$ are pairwise compatible.
\end{definition}

\begin{definition}
Let $A \in \mathrm{Ob}(\mathcal{C})$ and let $\dot{a} \in W^{(\mathcal{C},J)} (A)$.
$\dot{a}$ is \textit{atom type} if $\dot{a} = x^{(\mathrm{atom})}$ for some $x$.
$\dot{a}$ is \textit{set type} if $\dot{a} = x^{(\mathrm{set})}$ for some $x$.
When $\dot{a} = x^{(\mathrm{atom})}$ or  $\dot{a} = x^{(\mathrm{set})}$, we will also write $\dot{a}$ for $x$ if there is no confusion.
\end{definition}

\begin{definition}
We define the forcing relation $A \Vdash_{(\mathcal{C}, J)} \text{``}\, \phi (\dot{a}_0, \dot{a}_1, \dots , \dot{a}_{n-1} ) \,\text{''}$ for a formula $\phi (x_0, x_1, \dots , x_{n-1} )$ of $\mathcal{L}_{\mathrm{atom}}$, $A \in \mathrm{Ob}(\mathcal{C})$, and $\dot{a}_0, \dot{a}_1, \dots , \dot{a}_{n-1} \in W^{(\mathcal{C}, J)}(A)$ as follows:
\begin{itemize}
\item $A \Vdash_{(\mathcal{C}, J)} \text{``}\,\dot{a} \colon \mathrm{atom}\,\text{''}$ if and only if
	\begin{enumerate}
	\item $\emptyset \in J(A)$ or
	\item $\dot{a}$ is atom type.
	\end{enumerate}
\item $A \Vdash_{(\mathcal{C}, J)} \text{``}\, \dot{a} \colon \mathrm{set}\,\text{''}$ if and only if
	\begin{enumerate}
	\item $\emptyset \in J(A)$ or
	\item $\dot{a}$ is set type.
	\end{enumerate}
\item $A \Vdash_{(\mathcal{C}, J)} \text{``}\, \dot{a} \in \dot{b}\,\text{''}$ if and only if
	\begin{enumerate}
	\item $\emptyset \in J(A)$ or
	\item 
		\begin{enumerate}
		\item $\dot{b} \text{ is set type}$ and
		\item $\exists S \in J(A)\; \forall f \in S\; \exists \dot{x} \in W^{(\mathcal{C}, J)}(\dom f)$
			\begin{enumerate}
			\item $(\dot{x}, f) \in \dot{b}$ and
			\item $\dom f \Vdash_{(\mathcal{C}, J)} \text{``}\, \dot{a} \cdot f = \dot{x} \,\text{''}$.
			\end{enumerate}
		\end{enumerate}
	\end{enumerate}
\item $A \Vdash_{(\mathcal{C}, J)} \text{``}\, \dot{a} = \dot{b}\,\text{''}$ if and only if
	\begin{enumerate}
	\item $\emptyset \in J(A)$,
	\item 
		\begin{enumerate}
		\item $\dot{a}$ and $\dot{b}$ are atom type and
		\item $\exists S \in J(A)\; \forall f \in S\; (\dot{a} \cdot f = \dot{b} \cdot f)$, or
		\end{enumerate}
	\item 
		\begin{enumerate}
		\item $\dot{a}$ and $\dot{b}$ are set type and
		\item $\forall f \in \Hom_\mathcal{C} (\mathrm{any}, A)\; \forall \dot{x} \in W^{(\mathcal{C}, J)}(\dom f)$
			\begin{enumerate}
			\item $(\dot{x}, f) \in \dot{a} \rightarrow \dom f \Vdash_{(\mathcal{C}, J)} \text{``}\, \dot{x} \in \dot{b} \cdot f \,\text{''}$ and
			\item $(\dot{x}, f) \in \dot{b} \rightarrow \dom f \Vdash_{(\mathcal{C}, J)} \text{``}\,  \dot{x} \in \dot{a} \cdot f \,\text{''}$.
			\end{enumerate}
		\end{enumerate}
	\end{enumerate}
\item $A \Vdash_{(\mathcal{C},J) } \text{``}\, ( \phi \wedge \psi )( \dot{a}_0 , \dots , \dot{a}_{n-1} ) \,\text{''}$ if and only if
	\begin{enumerate}
	\item $A \Vdash_{(\mathcal{C},J) } \text{``}\,  \phi ( \dot{a}_0 , \dots , \dot{a}_{n-1} ) \,\text{''}$ and
	\item $A \Vdash_{(\mathcal{C},J) } \text{``}\,  \psi ( \dot{a}_0 , \dots , \dot{a}_{n-1} ) \,\text{''}$.
	\end{enumerate}
\item $A \Vdash_{(\mathcal{C},J) } \text{``}\, ( \phi \vee \psi )( \dot{a}_0 , \dots , \dot{a}_{n-1} ) \,\text{''}$ if and only if
	\begin{enumerate}
	\item[] $\exists S \in J(A)\; \forall f \in S$
		\begin{enumerate}
		\item $\dom f \Vdash_{(\mathcal{C},J) } \text{``}\,  \phi ( \dot{a}_0 \cdot f, \dots , \dot{a}_{n-1} \cdot f ) \,\text{''}$ or
		\item $\dom f \Vdash_{(\mathcal{C},J) } \text{``}\,  \psi ( \dot{a}_0 \cdot f, \dots , \dot{a}_{n-1} \cdot f ) \,\text{''}$.
		\end{enumerate}
	\end{enumerate}
\item $A \Vdash_{(\mathcal{C},J) } \text{``}\, ( \phi \rightarrow \psi )( \dot{a}_0 , \dots , \dot{a}_{n-1} ) \,\text{''}$ if and only if
	\begin{enumerate}
	\item[] $\forall f \in \Hom_\mathcal{C} (\mathrm{any}, A)$
		\begin{enumerate}
		\item $\dom f \Vdash_{(\mathcal{C}, J)} \text{``}\,  \phi ( \dot{a}_0 \cdot f, \dots , \dot{a}_{n-1} \cdot f) \,\text{''}$ implies
		\item $\dom f \Vdash_{(\mathcal{C}, J)} \text{``}\,  \psi ( \dot{a}_0 \cdot f, \dots , \dot{a}_{n-1} \cdot f) \,\text{''}$.
		\end{enumerate}
	\end{enumerate}
\item $A \Vdash_{(\mathcal{C},J) } \text{``}\,  \neg\phi ( \dot{a}_0 , \dots , \dot{a}_{n-1} ) \,\text{''}$ if and only if
	\begin{enumerate}
	\item[] $\forall f \in \Hom_\mathcal{C} (\mathrm{any}, A)$
		\begin{enumerate}
		\item $\dom f \Vdash_{(\mathcal{C}, J)} \text{``}\,  \phi ( \dot{a}_0 \cdot f, \dots , \dot{a}_{n-1} \cdot f) \,\text{''}$ implies
		\item $\emptyset \in J( \dom f ).$
		\end{enumerate}
	\end{enumerate}
\item $A \Vdash_{(\mathcal{C},J) } \text{``}\,  \forall x \phi ( x, \dot{a}_0 , \dots , \dot{a}_{n-1} ) \,\text{''}$ if and only if
	\begin{enumerate}
	\item[] $\forall f \in \Hom_\mathcal{C} (\mathrm{any}, A) \;\forall \dot{x} \in W^{(\mathcal{C},J)} (\dom f )$
		\begin{itemize}
		\item[] $\dom f \Vdash_{(\mathcal{C},J) } \text{``}\,  \phi ( \dot{x}, \dot{a}_0 \cdot f , \dots , \dot{a}_{n-1} \cdot f ) \,\text{''}$.
		\end{itemize}
	\end{enumerate}
\item $A \Vdash_{(\mathcal{C},J) } \text{``}\,  \exists x \phi ( x, \dot{a}_0 , \dots , \dot{a}_{n-1} ) \,\text{''}$ if and only if
	\begin{enumerate}
	\item[] $\exists S \in J(A) \;\forall f \in S \;\exists \dot{x} \in W^{(\mathcal{C},J)} (\dom f )$
		\begin{itemize}
		\item[] $\dom f \Vdash_{(\mathcal{C},J) } \text{``}\,  \phi ( \dot{x}, \dot{a}_0 \cdot f , \dots , \dot{a}_{n-1} \cdot f ) \,\text{''}$.
		\end{itemize}
	\end{enumerate}
\end{itemize}
\end{definition}

\subsection{Soundness}
\begin{proposition}\label{prop:forcingrel}
Let $\phi (x_0, \dots , x_{n-1} )$ be a formula of $\mathcal{L}_{\mathrm{atom}}$ and let $A \in \mathrm{Ob}(\mathcal{C})$. Let $\dot{a}_0, \dot{a}_1, \dots , \dot{a}_{n-1} \in W^{(\mathcal{C}, J)}(A)$.
\begin{enumerate}
\item If $A \Vdash_{(\mathcal{C}, J)} \text{``}\, \phi (\dot{a}_0 , \dots , \dot{a}_{n-1}) \text{''}$ holds,
then
\[
\dom f \Vdash_{(\mathcal{C}, J)} \text{``}\, \phi (\dot{a}_0 \cdot f , \dots , \dot{a}_{n-1} \cdot f) \text{''}
\]
holds for every $f \in \Hom_{\mathcal{C}} (\mathrm{any}, A)$.
\item If there exists $S \in J(A)$ for which
\[
\dom f \Vdash_{(\mathcal{C}, J)} \text{``}\, \phi (\dot{a}_0 \cdot f , \dots , \dot{a}_{n-1} \cdot f) \text{''}
\]
holds for every $f \in S$, then
$A \Vdash_{(\mathcal{C}, J)} \text{``}\, \phi (\dot{a}_0 , \dots , \dot{a}_{n-1}) \text{''}$ holds.
\end{enumerate}
\end{proposition}
\begin{proof}
By induction on $\phi (x_0, \dots , x_{n-1} )$.
\end{proof}

The following complete Heyting algebra $\Omega^{(\mathcal{C}, J)}(A)$ is convenient for describing propositions on the forcing relation.
\begin{definition}
Let $A \in \Ob (\mathcal{C})$.
A sieve $S$ on $A$ is $J$-\textit{closed} if for every $f \in \Hom_{\mathcal{C}}(\mathrm{any}, A)$, $f^* (S) \in J(\dom f)$ implies $f \in S$.
We define $\Omega^{(\mathcal{C}, J)}(A)$ to be the set of all $J$-closed sieves on $A$.
\end{definition}

\begin{proposition}
For every $A \in \Ob (\mathcal{C})$, the poset $\left( \Omega^{(\mathcal{C}, J)}(A), \subseteq \right)$ is a complete Heyting algebra in which the following properties hold:
\begin{enumerate}
\item $\bigwedge_{i \in I} S_i = \bigcap_{i \in I} S_i$,
\item $\bigvee_{i \in I} S_i = \left\{ f \in \Hom_{\mathcal{C}}(\mathrm{any}, A) \setmid f^* \left( \bigcup_{i \in I} S_i \right) \in J(\dom f) \right\}$,
\item $S_0 \rightarrow S_1 = \left\{ f \in \Hom_{\mathcal{C}}(\mathrm{any}, A) \setmid f^* (S_0) \subseteq f^* (S_1) \right\}$,
\item $1 = \Hom_{\mathcal{C}}(\mathrm{any}, A)$,
\item $0 = \left\{ f \in \Hom_{\mathcal{C}}(\mathrm{any}, A) \setmid \emptyset \in J(\dom f) \right\}$.
\end{enumerate}
\end{proposition}
\begin{proof}
Straightforward.
\end{proof}

\begin{definition}
Let  $\phi (x_0, \dots , x_{n-1} )$ be a formula of $\mathcal{L}_{\mathrm{atom}}$ and let $A \in \mathrm{Ob}(\mathcal{C})$. Let $\dot{a}_0, \dots , \dot{a}_{n-1} \in W^{(\mathcal{C}, J)}(A)$.
\begin{multline*}
\left\| \phi (\dot{a}_0 , \dots , \dot{a}_{n-1}) \right\|_A^{(\mathcal{C}, J)} := \{ f \in \Hom_\mathcal{C} (\mathrm{any}, A) \;|\\
\dom f \Vdash_{(\mathcal{C}, J)} \text{``}\,  \phi ( \dot{a}_0 \cdot f , \dots , \dot{a}_{n-1} \cdot f ) \,\text{''} \}.
\end{multline*}
\end{definition}

\begin{corollary}\label{cor:truthvalue}
Let $\phi (x_0, \dots , x_{n-1} )$ be a formula of $\mathcal{L}_{\mathrm{atom}}$ and let $A \in \mathrm{Ob}(\mathcal{C})$. Let $\dot{a}_0, \dot{a}_1, \dots , \dot{a}_{n-1} \in W^{(\mathcal{C}, J)}(A)$.
Then $\left\| \phi (\dot{a}_0 , \dots , \dot{a}_{n-1}) \right\|_A^{(\mathcal{C}, J)}$ is a $J$-closed sieve on $A$ i.e. $\left\| \phi (\dot{a}_0 , \dots , \dot{a}_{n-1}) \right\|_A^{(\mathcal{C}, J)} \in \Omega^{(\mathcal{C}, J)}(A)$.
\end{corollary}
\begin{proof}
Immediate from Proposition \ref{prop:forcingrel}.
\end{proof}

\begin{proposition}\label{prop:logico}
Let $\phi (x_0, \dots , x_{n-1} )$ and $\psi (x_0, \dots , x_{n-1} )$ be formulas of $\mathcal{L}_{\mathrm{atom}}$.
Let $A \in \mathrm{Ob}(\mathcal{C})$ and let $\dot{a}_0, \dots , \dot{a}_{n-1} \in W^{(\mathcal{C}, J)}(A)$.
Then in the complete Heyting algebra $\Omega^{(\mathcal{C}, J)}(A)$,
\begin{enumerate}
\item $\| ( \phi \vee \psi ) ( \dot{a}_0, \dots , \dot{a}_{n-1} ) \|_A^{(\mathcal{C},J)}\\
= \| \phi ( \dot{a}_0, \dots , \dot{a}_{n-1} ) \|_A^{(\mathcal{C},J)} \vee \| \psi ( \dot{a}_0, \dots , \dot{a}_{n-1} ) \|_A^{(\mathcal{C},J)}$,
\item $\| ( \phi \wedge \psi ) ( \dot{a}_0, \dots , \dot{a}_{n-1} ) \|_A^{(\mathcal{C},J)}\\
= \| \phi ( \dot{a}_0, \dots , \dot{a}_{n-1} ) \|_A^{(\mathcal{C},J)} \wedge \| \psi ( \dot{a}_0, \dots , \dot{a}_{n-1} ) \|_A^{(\mathcal{C},J)}$,
\item $\| ( \phi \rightarrow \psi ) ( \dot{a}_0, \dots , \dot{a}_{n-1} ) \|_A^{(\mathcal{C},J)}\\
= \| \phi ( \dot{a}_0, \dots , \dot{a}_{n-1} ) \|_A^{(\mathcal{C},J)} \rightarrow \| \psi ( \dot{a}_0, \dots , \dot{a}_{n-1} ) \|_A^{(\mathcal{C},J)}$,
\item $\| \neg \phi ( \dot{a}_0, \dots , \dot{a}_{n-1} ) \|_A^{(\mathcal{C},J)} = \neg \| \phi ( \dot{a}_0, \dots , \dot{a}_{n-1} ) \|_A^{(\mathcal{C},J)}$.
\end{enumerate}
\end{proposition}
\begin{proof}
Straightforward by the definition of the forcing relation.
\end{proof}

\begin{proposition}\label{prop:quo}
Let $\phi ( x, y_0, \dots , y_{n-1} )$ and $\psi (y_0, \dots , y_{n-1} )$ be formulas of $\mathcal{L}_{\mathrm{atom}}$. Let $A \in \mathrm{Ob}(\mathcal{C})$ and let $\dot{a}_0, \dots , \dot{a}_{n-1}, \dot{b} \in W^{(\mathcal{C}, J)}(A)$.
\begin{enumerate}
\item $\| \forall x \phi (x, \dot{a}_0 , \dots , \dot{a}_{n-1}) \|_A^{(\mathcal{C}, J)} \le \| \phi (\dot{b}, \dot{a}_0 , \dots , \dot{a}_{n-1}) \|_A^{(\mathcal{C}, J)}$,
\item $\| \phi (\dot{b}, \dot{a}_0 , \dots , \dot{a}_{n-1}) \|_A^{(\mathcal{C}, J)} \le \| \exists x \phi (x, \dot{a}_0 , \dots , \dot{a}_{n-1}) \|_A^{(\mathcal{C}, J)}$,
\item $\| \forall x ( \psi ( \dot{a}_0, \dots , \dot{a}_{n-1} ) \rightarrow \phi (x, \dot{a}_0 , \dots , \dot{a}_{n-1}) ) \|_A^{(\mathcal{C}, J)}\\
\le \| \psi ( \dot{a}_0, \dots , \dot{a}_{n-1} ) \rightarrow \forall x \phi (x, \dot{a}_0 , \dots , \dot{a}_{n-1}) \|_A^{(\mathcal{C}, J)}$,
\item $\| \forall x ( \phi (x, \dot{a}_0, \dots , \dot{a}_{n-1} ) \rightarrow \psi ( \dot{a}_0 , \dots , \dot{a}_{n-1}) ) \|_A^{(\mathcal{C}, J)}\\
\le \| (\exists x \phi (x, \dot{a}_0, \dots , \dot{a}_{n-1} )) \rightarrow \psi ( \dot{a}_0 , \dots , \dot{a}_{n-1}) \|_A^{(\mathcal{C}, J)}$.
\end{enumerate}
\end{proposition}
\begin{proof}
Straightforward by the definition of the forcing relation.
\end{proof}

\begin{proposition}\label{prop:gen}
Let $\phi ( x, y_0 , \dots , y_{n-1})$ be a formula of $\mathcal{L}_{\mathrm{atom}}$.
If
\[
\| \phi (\dot{a}, \dot{b}_0, \dots , \dot{b}_{n-1}) \|_A^{(\mathcal{C}, J)} = 1
\]
holds for every $A \in \mathrm{Ob}(\mathcal{C})$ and every $\dot{a}, \dot{b}_0, \dots , \dot{b}_{n-1} \in W^{(\mathcal{C}, J)}(A)$,
then
\[
\| \forall x \phi (x, \dot{b}_0, \dots , \dot{b}_{n-1}) \|_A^{(\mathcal{C}, J)} = 1
\]
holds for every $A \in \mathrm{Ob}(\mathcal{C})$ and every $\dot{b}_0, \dots , \dot{b}_{n-1} \in W^{(\mathcal{C}, J)}(A)$.
\end{proposition}
\begin{proof}
Straightforward by the definition of the forcing relation.
\end{proof}

\begin{proposition}\label{prop:equ}
Let $A \in \mathrm{Ob}(\mathcal{C})$ and let $\dot{a}, \dot{b}, \dot{c} \in W^{(\mathcal{C},J)}(A)$.
\begin{enumerate}
\item $\| \dot{a} = \dot{a} \|_A^{(\mathcal{C}, J)} = 1$,\label{enu:refl}
\item $\| \dot{a} = \dot{b} \|_A^{(\mathcal{C}, J)} \le \| \dot{b} = \dot{a} \|_A^{(\mathcal{C}, J)} $,\label{enu:sym}
\item $\| \dot{a} = \dot{b} \|_A^{(\mathcal{C}, J)} \wedge \| \dot{b} = \dot{c} \|_A^{(\mathcal{C}, J)} \le \| \dot{a} = \dot{c} \|_A^{(\mathcal{C}, J)}$,\label{enu:trans}
\item $\| \dot{a} \in \dot{b} \|_A^{(\mathcal{C}, J)} \wedge \| \dot{a} = \dot{c} \|_A^{(\mathcal{C}, J)} \le \| \dot{c} \in \dot{b} \|_A^{(\mathcal{C}, J)}$,\label{enu:submem1}
\item $\| \dot{a} \in \dot{b} \|_A^{(\mathcal{C}, J)} \wedge \| \dot{b} = \dot{c} \|_A^{(\mathcal{C}, J)} \le \| \dot{a} \in \dot{c} \|_A^{(\mathcal{C}, J)}$,\label{enu:submem2}
\item $\| \dot{a} \colon \mathrm{atom} \|_A^{(\mathcal{C}, J)} \wedge \| \dot{a} = \dot{b} \|_A^{(\mathcal{C}, J)} \le \| \dot{b} \colon \mathrm{atom} \|_A^{(\mathcal{C}, J)}$,\label{enu:subatom}
\item $\| \dot{a} \colon \mathrm{set} \|_A^{(\mathcal{C}, J)} \wedge \| \dot{a} = \dot{b} \|_A^{(\mathcal{C}, J)} \le \| \dot{b} \colon \mathrm{set} \|_A^{(\mathcal{C}, J)}$.\label{enu:subset}
\end{enumerate}
\end{proposition}
\begin{proof}
(1): By induction on $\dot{a}$.

(2), (6), (7): Straightforward by the definition of the forcing relation.

(3), (4), (5): By simultaneous induction on $\dot{a}, \dot{b}, \dot{c}$.
\end{proof}

\begin{theorem}[Soundness]
Let $\phi (x_0 , \dots , x_{n-1})$ be a formula of $\mathcal{L}_\mathrm{atom}$.
If $\phi$ is provable in intuitionistic first-order logic with equality,
then
\[
A \Vdash_{(\mathcal{C}, J)} \text{``}\, \phi (\dot{a}_0 , \dots , \dot{a}_{n-1}) \,\text{''}
\]
holds for every $A \in \mathrm{Ob}(\mathcal{C})$ and every $\dot{a}_0 , \dots , \dot{a}_{n-1} \in W^{(\mathcal{C}, J)} (A)$.
\end{theorem}
\begin{proof}
It is sufficient to show that $\| \phi (\dot{a}_0, \dots , \dot{a}_{n-1}) \|_A^{(\mathcal{C}, J)} = 1$ for all $A \in \mathrm{Ob}(\mathcal{C})$ and all $\dot{a}_0 , \dots , \dot{a}_{n-1} \in W^{(\mathcal{C}, J)} (A)$, but it is straightforward by Propositions \ref{prop:logico}, \ref{prop:quo}, \ref{prop:gen}, and \ref{prop:equ}.
\end{proof}

\subsection{Check operator}
\begin{definition}[Check operator]
For a set $x$ and an object $A \in \mathrm{Ob}(\mathcal{C})$, we recursively define $\Check{x}^A$ (or $(x)\Check{}^{A}$) $\in W^{(\mathcal{C}, J)} (A)$ by
\[
\Check{x}^A = \left\{ (\Check{y}^{\dom f}, f) \setmid y \in x,\; f \in \Hom_\mathcal{C} (\mathrm{any}, A) \right\} .
\]
\end{definition}

\begin{proposition}
Let $x$ be a set and let $A \in \mathrm{Ob}(\mathcal{C})$. Then $\Check{x}^A \cdot f = \Check{x}^{\dom f}$ for all $f \in \Hom_\mathcal{C} (\mathrm{any}, A)$.
\end{proposition}
\begin{proof}
Straightforward.
\end{proof}

\begin{theorem}
Let $\phi (x_0, \dots , x_{n-1})$ be a $\Delta_0$-formula of $\mathcal{L}_\in$ and let $A \in \mathrm{Ob}(\mathcal{C})$. Let $a_0, \dots , a_{n-1}$ be sets. Then
\[
\| \phi (\Check{a}_0^A , \dots ,\Check{a}_{n-1}^A ) \|_A^{(\mathcal{C},J)} =
\begin{cases}
1 & \text{if } \phi (a_0, \dots , a_{n-1}) \text{ holds } (\text{in } V), \\
0 & \text{otherwise.}
\end{cases}
\]
\end{theorem}
\begin{proof}
By induction on $\phi (x_0, \dots , x_{n-1})$.
\end{proof}

\subsection{Maximum principle}
\begin{definition}
Let $A \in \mathrm{Ob}(\mathcal{C})$ and let $S$ be a sieve on $A$.
A function $F$ on $S$ is called a \textit{matching function for} $S$ if the following conditions hold:
\begin{enumerate}
\item $F(f)$ is a nonempty subset of $W^{(\mathcal{C}, J)}(\dom f)$ for every $f \in S$,
\item For every $f \in S$ and every $g \in \Hom_\mathcal{C} (\mathrm{any}, \dom f)$, if $\dot{a} \in F(f)$ and $\dot{b} \in F(f \circ g)$, then $\dom g \Vdash_{(\mathcal{C}, J)} \text{``}\,  \dot{a} \cdot g = \dot{b} \,\text{''}$.
\end{enumerate}
\end{definition}

\begin{definition}
Let $A \in \mathrm{Ob}(\mathcal{C})$ and let $S$ be a sieve on $A$. Let $F$ be a matching function for $S$.
We assume that all elements of $F(f)$ are set type for every $f \in S$. 
Then we define the \textit{amalgamation of} $F$ by
\[
\ama F := \left\{ (\dot{x}, f \circ g) \setmid f \in S,\; \exists \dot{a} \in F(f) \left( (\dot{x}, g) \in \dot{a} \right) \right\} \in W^{(\mathcal{C}, J)}(A).
\]
\end{definition}

\begin{proposition}\label{prop:ama}
Let $A \in \mathrm{Ob}(\mathcal{C})$ and let $S$ be a sieve on $A$. Let $F$ be a matching function for $S$.
We assume that all elements of $F(f)$ are set type for every $f \in S$. 
Let $\dot{a} := \ama F$.
Then $\dom f \Vdash_{(\mathcal{C}, J)} \text{``}\,  \dot{a} \cdot f = \dot{b} \,\text{''}$ for all $f \in S$ and all $\dot{b} \in F(f)$.
\end{proposition}
\begin{proof}
Straightforward by the definition of the forcing relation.
\end{proof}

\begin{theorem}[Maximum principle]\label{thm:max}
Let $\phi (x , y_0 , \dots , y_{n-1})$ be a formula of $\mathcal{L}_{atom}$ and let $A \in \mathrm{Ob}(\mathcal{C})$.
Let $\dot{a}_0 , \dots , \dot{a}_{n-1} \in W^{(\mathcal{C}, J)}(A)$. \\
If $A \Vdash_{(\mathcal{C}, J)} \text{``}\,  \exists ! x \colon \mathrm{set}, \phi (x, \dot{a}_0 , \dots , \dot{a}_{n-1}) \,\text{''}$,
then there exists $\dot{x} \in W^{(\mathcal{C}, J)}(A)$ for which $A \Vdash_{(\mathcal{C}, J)} \text{``}\, \dot{x} \colon \mathrm{set} \wedge \phi (\dot{x}, \dot{a}_0 , \dots , \dot{a}_{n-1}) \,\text{''}$.
\end{theorem}
\begin{proof}
Straightforward by Proposition \ref{prop:ama}.
\end{proof}

\subsection{IZFA}
\begin{definition}
\textit{Intuitionistic Zermelo-Fraenkel set theory with atoms} (or IZFA) is the theory in $\mathcal{L}_{\mathrm{atom}}$ based on the following axioms:
\begin{enumerate}
\item Set existence
\[
\exists x\; (x \colon \mathrm{set}).
\]
\item Extensionality
\[
\forall x \colon \mathrm{set} \; \forall y \colon \mathrm{set}\; (\forall z\; (z \in x \leftrightarrow z \in y) \rightarrow x=y).
\]
\item Separation
\[
\forall u \colon \mathrm{set}\; \exists v \colon \mathrm{set}\; \forall x\; (x \in v \leftrightarrow x \in u \wedge \phi (x)),
\]
where $v$ is not free in the formula $\phi (x)$ of $\mathcal{L}_{\mathrm{atom}}$.
\item Collection
\[
\forall u \colon \mathrm{set}\; ( \forall x \in u\; \exists y \; \phi (x, y) \rightarrow \exists v \colon \mathrm{set}\; \forall x \in u \; \exists y \in v\; \phi (x, y) ),
\]
where $v$ is not free in the formula $\phi (x, y)$ of $\mathcal{L}_{\mathrm{atom}}$.
\item Pairing \label{axiom:pairing}
\[
\forall x\; \forall y\; \exists z\; \forall w\; (w \in z \leftrightarrow w=x \vee w=y).
\]
\item Union
\[
\forall u \colon \mathrm{set}\; \exists v \colon \mathrm{set}\; \forall x\; ( x \in v \leftrightarrow \exists y \in u\; (x \in y) ).
\]
\item Power set
\[
\forall u \colon \mathrm{set}\; \exists v \colon \mathrm{set}\; \forall x\; ( x \in v \leftrightarrow \forall y \in x\; (y \in u) ).
\]
\item Infinity
\[
\exists u\; ( \emptyset \in u \wedge \forall x \in u\; ( x \cup \{ x \} \in u ) ).
\]
\item $\in$-induction
\[
\forall x \; ( \forall y \in x\; \phi (y) \rightarrow \phi (x) ) \rightarrow \forall x\; \phi (x),
\]
where $y$ is not free in the formula $\phi (x)$ of $\mathcal{L}_{\mathrm{atom}}$.
\item Atom
\[
\forall x \colon \mathrm{atom}\; \forall y\; (y \not\in x),
\]
\[
\forall x \; (x \colon \mathrm{atom} \vee x \colon \mathrm{set}),
\]
\[
\forall x \; \neg (x \colon \mathrm{atom} \wedge x \colon \mathrm{set}).
\]
\end{enumerate}
\end{definition}

\begin{definition}
For $A \in \Ob (\mathcal{C})$ and $\dot{a}, \dot{b} \in W^{(\mathcal{C}, J)}(A)$, we define the \textit{unordered pair} and the \textit{ordered pair} of $\dot{a}, \dot{b}$ in $W^{(\mathcal{C}, J)}(A)$ as follows:
\begin{itemize}
\item $\{ \dot{a}, \dot{b} \}_{A} := \setbegin ( \dot{a} \cdot f, f ) \setmid f \in \Hom_{\mathcal{C}}(\mathrm{any}, A) \setend \cup \{\, ( \dot{b} \cdot f, f ) \, | \, f \in \Hom_{\mathcal{C}}(\mathrm{any}, A) \,\}$,
\item $(\dot{a}, \dot{b})_A := \{ \{ \dot{a}, \dot{a} \}_{A} , \{ \dot{a}, \dot{b} \}_{A} \}_{A}$.
\end{itemize}
\end{definition}

\begin{theorem}
For every axioms $\phi$ of $\mathrm{IZFA}$ and every $A \in \mathrm{Ob}(\mathcal{C})$,
\[
A \Vdash_{(\mathcal{C}, J)} \text{``}\, \phi \,\text{''}
\]
holds.
\end{theorem}
\begin{proof}
Easy. For example, $\dot{z} := \{ \dot{x}, \dot{y} \}_A$ is a witness for (\ref{axiom:pairing}) Pairing.
\end{proof}

\section{Toposes from forcing}\label{sec:toposes}
Let $(\mathcal{C}, J)$ be a site.
In this section, we define the category $\mathrm{Set}^{(\mathcal{C}, J)}$ of `sets in the forcing model $(W^{(\mathcal{C}, J)}, \Vdash_{(\mathcal{C}, J)})$' and prove the main theorem that the categories $\Sh (\mathcal{C}, J)$ and $\mathrm{Set}^{(\mathcal{C}, J)}$ are equivalent by constructing a fully faithful and essentially surjective functor $L \colon \mathrm{Set}^{(\mathcal{C}, J)} \rightarrow \mathrm{Sh} (\mathcal{C}, J)$.
Henceforth, for each $A \in \Ob (\mathcal{C})$, we identify elements $\dot{a}, \dot{b}$ of $W^{(\mathcal{C}, J)} (A)$ when $A\Vdash_{(\mathcal{C}, J)} \text{``}\,  \dot{a} = \dot{b} \,\text{''}$.

\subsection{Category $\mathrm{Set}^{(\mathcal{C}, J)}$ of `sets in $(W^{(\mathcal{C}, J)}, \Vdash_{(\mathcal{C}, J)})$'}
\begin{definition}\text{}
\begin{itemize}
\item A $(\mathcal{C}, J)$-\textit{sequence} is a sequence $(\dot{a}_A)_{A \in \Ob (\mathcal{C})}$ of which each $\dot{a}_A$ is an element of $ W^{(\mathcal{C}, J)}(A)$.
\item A $(\mathcal{C}, J)$-sequence $(\dot{a}_A)_{A \in \Ob (\mathcal{C})}$ is called \textit{stable} if
\[
\dom f \Vdash_{(\mathcal{C}, J)} \text{``}\, \dot{a}_{\cod f} \cdot f = \dot{a}_{\dom f} \,\text{''}
\]
holds for every $f \in \Arr (\mathcal{C})$.
\item A $(\mathcal{C}, J)$-\textit{set} is a stable $(\mathcal{C}, J)$-sequence $(\dot{a}_A)_{A \in \Ob (\mathcal{C})}$ for which
\[
A\Vdash_{(\mathcal{C}, J)} \text{``}\,  \dot{a}_A \colon \mathrm{set} \,\text{''}
\]
holds for every $A \in \Ob (\mathcal{C})$.
\end{itemize}
\end{definition}

\begin{definition}
We define a category $\mathrm{Set}^{(\mathcal{C}, J)}$ as follows:
\begin{itemize}
\item the objects of $\mathrm{Set}^{(\mathcal{C}, J)}$ are the $(\mathcal{C}, J)$-sets,
\item the arrows of $\mathrm{Set}^{(\mathcal{C}, J)}$ from $(\dot{a}_A)_{A \in \Ob (\mathcal{C})}$ to $(\dot{b}_A)_{A \in \Ob (\mathcal{C})}$ are those $(\mathcal{C}, J)$-sets $(\dot{f}_A)_{A \in \Ob (\mathcal{C})}$ for which $A \Vdash_{(\mathcal{C}, J)} \text{``}\, \dot{f}_A \text{ is a function from } \dot{a}_A \text{ to } \dot{b}_A \,\text{''}$ for every $A \in \Ob (\mathcal{C})$,
\item the composition of two arrows $(\dot{f}_A)_{A \in \Ob (\mathcal{C})}$ and $(\dot{g}_A)_{A \in \Ob (\mathcal{C})}$ of $\mathrm{Set}^{(\mathcal{C}, J)}$ is the unique arrow $(\dot{h}_A)_{A \in \Ob (\mathcal{C})}$ for which $A \Vdash_{(\mathcal{C},J)} \text{``}\, \dot{f}_A \circ \dot{g}_A = \dot{h}_A \,\text{''}$ for every $A \in \Ob (\mathcal{C})$. (Such $\dot{h}_A$ exists by Theorem \ref{thm:max}.)
\end{itemize}
\end{definition}

\subsection{Functor $L \colon \mathrm{Set}^{(\mathcal{C}, J)} \rightarrow \mathrm{Sh} (\mathcal{C}, J)$}
\begin{definition}
Let $\mathbf{a} = (\dot{a}_A)_{A \in \Ob (\mathcal{C})}$ be a $(\mathcal{C}, J)$-set.
We define a presheaf $L_{\mathbf{a}}^{\mathrm{pre}}$ on $\mathcal{C}$ as follows:
\begin{itemize}
\item $L_{\mathbf{a}}^{\mathrm{pre}} (A) := \setbegin \dot{c} \in W^{(\mathcal{C}, J)} (A) \setmid A \Vdash_{(\mathcal{C}, J)} \text{``}\, \dot{c} \in \dot{a}_A \,\text{''} \setend$ for $A \in \Ob (\mathcal{C})$,
\item for $f \in \Hom_{\mathcal{C}} (A, B)$, a function $L_{\mathbf{a}}^{\mathrm{pre}} (f) \colon L_{\mathbf{a}}^{\mathrm{pre}} (B) \to L_{\mathbf{a}}^{\mathrm{pre}} (A)$ is defined by $L_{\mathbf{a}}^{\mathrm{pre}} (f) (\dot{c}) = \dot{c} \cdot f$.
\end{itemize}
Let $L_{\mathbf{a}}$ be the sheafification of $L_{\mathbf{a}}^{\mathrm{pre}}$ and let $i_{\mathbf{a}} = ( i_{\mathbf{a}, A} )_{A \in \Ob (\mathcal{C})} \colon L_{\mathbf{a}}^{\mathrm{pre}} \to L_{\mathbf{a}}$ be its canonical map.
Since $L_{\mathbf{a}}^{\mathrm{pre}}$ is a separated presheaf, $i_{\mathbf{a}}$ is a monomorphism.
\end{definition}

\begin{definition}
Let $\mathbf{f} = ( \dot{f}_A )_{A \in \Ob (\mathcal{C})} \in \Hom_{\mathrm{Set}^{(\mathcal{C}, J)}} (\mathbf{a}, \mathbf{b})$.
We define a natural transformation $L_{\mathbf{f}}^{\mathrm{pre}} = \left( L_{\mathbf{f}, A}^{\mathrm{pre}} \right)_{A \in \Ob (\mathcal{C})} \colon L_{\mathbf{a}}^{\mathrm{pre}} \to L_{\mathbf{b}}^{\mathrm{pre}}$ by:
\[
L_{\mathbf{f}, A}^{\mathrm{pre}} (\dot{c}) = \dot{d} \ \ \text{if and only if}\ \ A \Vdash_{(\mathcal{C}, J)} \text{``}\, \dot{f}_A (\dot{c}) = \dot{d} \,\text{''.}
\]
Let $L_{\mathbf{f}}$ be the unique natural transformation $\sigma \colon L_{\mathbf{a}} \to L_{\mathbf{b}}$ for which the following diagram commutes:
\[
\begin{CD}
L_{\mathbf{a}}^{\mathrm{pre}} @> i_{\mathbf{a}} >> L_{\mathbf{a}} \\
@V L_{\mathbf{f}}^{\mathrm{pre}} VV @VV \sigma V \\
L_{\mathbf{b}}^{\mathrm{pre}} @>> i_{\mathbf{b}} > L_{\mathbf{b}}
\end{CD}
\]
\end{definition}

\begin{definition}
We define a functor $L \colon \mathrm{Set}^{(\mathcal{C}, J)} \rightarrow \mathrm{Sh} (\mathcal{C}, J)$ as follows:
\begin{itemize}
\item $L (\mathbf{a}) := L_{\mathbf{a}}$ for $\mathbf{a} \in \Ob (\mathrm{Set}^{(\mathcal{C}, J)})$,
\item $L (\mathbf{f}) := L_{\mathbf{f}}$ for $\mathbf{f} \in \Hom_{\mathrm{Set}^{(\mathcal{C}, J)}} (\mathbf{a}, \mathbf{b})$.
\end{itemize}
\end{definition}

\subsection{Representation of sheaves on $(\mathcal{C}, J)$ by $(\mathcal{C}, J)$-sets}
Before proving the main theorem, we construct the following $(\mathcal{C}, J)$-set $\bigl( \overline{F}^A \bigr)_{A \in \Ob (\mathcal{C})}$ for each sheaf $F$ on $(\mathcal{C}, J)$, which is used for showing that the functor $L$ is essentially surjective.
\begin{definition}
Let $F$ be a sheaf on $(\mathcal{C}, J)$ and $A \in \Ob (\mathcal{C})$.
For each $a \in F(A)$, we define
\begin{multline*}
\overline{a}^{F, A} := \setbegin \left( \left( \Check{x}^{\dom g} , f^{(\mathrm{atom})} \right)_{\dom g},\; g \right) \setmid f \in \Arr (\mathcal {C}),\; x \in F(\cod f),\right.\\
g \in \Hom_\mathcal{C} (\dom f, A),\; F(f)(x) = F(g)(a) \,\Biggr\} \in W^{(\mathcal{C}, J)}(A)
\end{multline*}
and let
$\overline{F}^A := \setbegin \left( \overline{a}^{F, \dom f}, f \right) \setmid f \in \Hom_{\mathcal{C}} (\mathrm{any}, A),\, a \in F ( \dom f ) \setend \in W^{(\mathcal{C}, J)}(A)$.
\end{definition}

Then these elements $\overline{a}^{F, A}$, $\overline{F}^A$ of $W^{(\mathcal{C}, J)} (A)$ represent the behavior of a sheaf $F$ well as follows:
\begin{proposition} \label{prop:enc}
Let $F$ be a sheaf on $(\mathcal{C}, J)$ and let $A \in \Ob (\mathcal{C})$.
\begin{enumerate}
\item \label{prop:enc:set}$\left( \overline{F}^A \right)_{A \in \Ob (\mathcal{C})}$ is a $(\mathcal{C}, J)$-set.
\item \label{prop:enc:ele}$A\Vdash_{(\mathcal{C}, J)} \text{``}\ \overline{\mathstrut a}^{F, A} \in \overline{\mathstrut F}^A \text{''}$ for all $a \in F(A)$.
\item \label{prop:enc:pres}$\overline{a}^{F, A} \cdot h = \overline{F(h)(a)}^{F, \dom h}$ for all $a \in F(A)$ and all $h \in \Hom_\mathcal{C} (\mathrm{any}, A)$.
\item \label{prop:enc:inj}For $a, b \in F(A)$, if $A\Vdash_{(\mathcal{C}, J)} \text{``}\ \overline{\mathstrut a}^{F, A} = \overline{\mathstrut b}^{F, A} \text{''}$,
then $a = b$.
\item \label{prop:enc:surj}For $\dot{x} \in W^{(\mathcal{C}, J)}(A)$, if $A\Vdash_{(\mathcal{C}, J)} \text{``}\ \dot{x} \in \overline{F}^A \text{''}$,
then there exists $a \in F(A)$ for which $A\Vdash_{(\mathcal{C}, J)} \text{``}\ \dot{x} = \overline{a}^{F, A} \text{''}$.
\end{enumerate}
\end{proposition}
\begin{proof}
(1), (2), (3): Immediate.

(4): Straightforward since $A\Vdash_{(\mathcal{C}, J)} \text{``}\, \left(\Check{a}^A, 1_A^{(\mathrm{atom})} \right) \in \overline{a}^{F, A} \text{''}$.

(5): Straightforward by (3) and (4) since $F$ is a sheaf on $(\mathcal{C}, J)$.
\end{proof}

\subsection{Main theorem}
Now we will prove the main theorem:
\begin{theorem}\label{thm:equiv}
The functor $L \colon \mathrm{Set}^{(\mathcal{C}, J)} \rightarrow \mathrm{Sh} (\mathcal{C}, J)$ is fully faithful and essentially surjective.
Thus, $\mathrm{Set}^{(\mathcal{C}, J)}$ and $\mathrm{Sh} (\mathcal{C}, J)$ are equivalent.
\end{theorem}
\begin{proof}

[Fullness of $L$]:
Let $\mathbf{a}$ and $\mathbf{b}$ be $(\mathcal{C}, J)$-sets and let $\sigma = ( \sigma_A )_{A \in \Ob (\mathcal{C})}$ be a natural transformation from $L_{\mathbf{a}}$ to $L_{\mathbf{b}}$.
For each $A \in \Ob (\mathcal{C}, J)$ we define
\begin{multline*}
\dot{f}_A := \setbegin \left( (\dot{c}, \dot{d})_A \cdot g,\, g \right) \setmid \dot{c} \in L_{\mathbf{a}}^{\mathrm{pre}} (A),\, \dot{d} \in L_{\mathbf{b}}^{\mathrm{pre}} (A),\right.\\
\left. g \in \Hom_{\mathcal{C}} (\mathrm{any}, A),\, \sigma_A \circ i_{\mathbf{a}, A} (\dot{c}) = i_{\mathbf{b}, A} (\dot{d}) \setend \in W^{(\mathcal{C}, J)} (A).
\end{multline*}
and let $\mathbf{f} := ( \dot{f}_A )_{A \in \Ob (\mathcal{C})}$.
Then we can prove easily that $\mathbf{f}$ is an arrow of $\mathrm{Set}^{(\mathcal{C}, J)}$ from $\mathbf{a}$ to $\mathbf{b}$.
By the definitions of $L_{\mathbf{f}}$, $L_{\mathbf{f}}^{\mathrm{pre}}$ and $\mathbf{f}$, it holds that
\[
L_{\mathbf{f}} \circ i_{\mathbf{a}} = i_{\mathbf{b}} \circ L_{\mathbf{f}}^{\mathrm{pre}} = \sigma \circ i_{\mathbf{a}}.
\]
Hence, $L_{\mathbf{f}} = \sigma$ by the universal property for the canonical map $i_{\mathbf{a}} \colon L_{\mathbf{a}}^{\mathrm{pre}} \to L_{\mathbf{a}}$.  

[Faithfulness of $L$]:
Let $\mathbf{f}$ and $\mathbf{g}$ be arrows of $\mathrm{Set}^{(\mathcal{C}, J)}$ from $\mathbf{a}$ to $\mathbf{b}$ for which $L_{\mathbf{f}} = L_{\mathbf{g}}$ holds.
Then
\[
i_{\mathbf{b}} \circ L_{\mathbf{f}}^{\mathrm{pre}} = L_{\mathbf{f}} \circ i_{\mathbf{a}} = L_{\mathbf{g}} \circ i_{\mathbf{a}} = i_{\mathbf{b}} \circ L_{\mathbf{g}}^{\mathrm{pre}}.
\]
Since $i_{\mathbf{b}}$ is a monomorphism, $L_{\mathbf{f}}^{\mathrm{pre}} = L_{\mathbf{g}}^{\mathrm{pre}}$.
Hence, by the definition of $L_{\mathbf{f}}^{\mathrm{pre}}$ and $L_{\mathbf{g}}^{\mathrm{pre}}$, $\mathbf{f} = \mathbf{g}$ holds.

[Essential surjectivity of $L$]:
Fix a sheaf $F$ on $(\mathcal{C}, J)$.
Let $\widetilde{F} := \left( \overline{F}^A \right)_{A \in \Ob (\mathcal{C})}$.
By Proposition \ref{prop:enc} (\ref{prop:enc:set}) and (\ref{prop:enc:ele}),
we can define a function $\sigma_A \colon F (A) \to L_{\widetilde{F}}^{\mathrm{pre}} (A)$ for each $A \in \Ob (\mathcal{C})$ by
\[
\sigma_A (a) = \overline{a}^{F, A}.
\]
Then $\sigma := ( \sigma_A )_{A \in \Ob (\mathcal{C})}$ has the following properties:
\begin{itemize}
\item $\sigma$ is a natural transformation from $F$ to $L_{\widetilde{F}}^{\mathrm{pre}}$ by Proposition \ref{prop:enc} (\ref{prop:enc:pres}),
\item each $\sigma_A \colon F (A) \to L_{\widetilde{F}}^{\mathrm{pre}} (A)$ is injective by Proposition \ref{prop:enc} (\ref{prop:enc:inj}),
\item it is also surjective by Proposition \ref{prop:enc} (\ref{prop:enc:surj}).
\end{itemize}
Hence, $\sigma$ is a natural isomorphism from $F$ to $L_{\widetilde{F}}^{\mathrm{pre}}$.
Since $F$ is a sheaf, $L_{\widetilde{F}}^{\mathrm{pre}}$ is also a sheaf, which is isomorphic to its sheafification $L_{\widetilde{F}}$.
Therefore, $F$ is isomorphic to $L_{\widetilde{F}}$. 
\end{proof}

\nocite{*}

\end{document}